\documentclass[11pt]{amsart}

\usepackage{amsmath,amssymb,epsf}
\usepackage{psfig}
\usepackage{amsfonts}
\usepackage{times}

\input xy
\xyoption{all}

\newtheorem{thm}{Theorem}[section]
\newtheorem{lem}[thm]{Lemma}
\newtheorem{cor}[thm]{Corollary}
\newtheorem{prop}[thm]{Proposition}
\newtheorem{exe}[thm]{Example}

\theoremstyle{definition}
\newtheorem{defn}[thm]{Definition}

\newcommand{\ZZ}{\mathbb{Z}}

\newcommand{\RR}{\mathbb{R}}
\newcommand{\CC}{\mathbb{C}}

\newcommand{\KK}[1]{\mathbb{K}_{{#1}}}

\newcommand{\MC}{\mathcal{C}}
\newcommand{\MF}{\mathcal{F}}
\newcommand{\MH}{\mathcal{H}}
\newcommand{\MV}{\mathcal{V}}

\newcommand{\suppfields}[3]{\overline{\MF}_{#1}^{#2,#3}}

\newcommand{\dprime}{{\prime\prime}}

\newcommand{\leg}[2]{(\frac{#1}{#2})}

\newcommand{\cohomol}[1]{H^1(#1;A)}

\newcommand{\ra}{\rightarrow}
\newcommand{\lra}{\longrightarrow}

\newcommand{\op}[1]{\overline{#1}}
\newcommand{\MFr}[3]{\MF_{#1}^{{#2},{#3}}}
\newcommand{\MFrr}[4]{\MF_{#1}^{{#2},{#3},{#4}}}
\newcommand{\MFLpq}{\MF_{L(p,q)}}
\newcommand{\meas}[1]{\mu_{\!{#1}}^{\empty}}
\newcommand{\measalpha}{\mu_{\!W}^{\alpha_0,\alpha_1}}
\newcommand{\measr}[3]{\mu_{\!{#1}}^{#2,#3}}
\newcommand{\suppmeas}[3]{\overline{\mu}_{\!{#1}}^{#2,#3}}
\newcommand{\dmeas}[1]{d\!\meas {#1}}
\newcommand{\dmeasalpha}{d\!\mu_{\!W}^{\alpha_0,\alpha_1}}

\newcommand{\dsuppmeas}[3]{d\!\overline{\mu}_{\!{#1}}^{#2,#3}}
\newcommand{\sigp}{\sigma^\prime}
\newcommand{\gpzero}{\gamma^\prime_0}
\newcommand{\gpone}{\gamma^\prime_1}
\newcommand{\slant}{\setminus}

\newcommand{\qb}{\bar{q}}

\DeclareMathOperator{\Hom}{Hom}
\DeclareMathOperator{\Ker}{Ker}
\DeclareMathOperator{\Id}{Id}

\DeclareMathOperator{\Map}{Map}

\title{Abelian homotopy Dijkgraaf--Witten theory}
\author{S.~K.~ Hansen}
\address{Department of Mathematics \\Kansas State University\\
  Manhattan \\ KS 66506\\ USA}
\author{J.~K.~ Slingerland}
\address{Microsoft Research\\1 Microsoft Way 113/2033\\
  Redmond \\ WA 98052\\ USA}
\author{P.~R.~Turner}
\address{School of Mathematical and Computer Sciences \\
  Heriot-Watt University\\ Edinburgh EH14 4AS \\Scotland}

\begin{document}


\begin{abstract}
We construct a version of Dijkgraaf--Witten theory based on a compact
abelian Lie group within the formalism of Turaev's homotopy quantum
field theory. As an application we show that the 2+1--dimensional
theory based on $U(1)$ classifies lens spaces up to homotopy type.
\end{abstract}

\maketitle


\section{Introduction}

The central topic of this paper is Dijkgraaf--Witten (DW)
invariants of closed, oriented $n+1$--manifolds based on
a compact abelian gauge group $A$. These may be defined as
follows.

The ``space of fields'' on an $n+1$--manifold $W$ is taken
to be the moduli space $\MF_W$ of isomorphism classes of
$A$--bundles with flat connection. Since $A$ is abelian
there are identifications
\[
\MF_W \cong \Hom(\pi_1(W),A)/\mbox{conj} \cong
\Hom(\pi_1(W),A) \cong \Hom(H_1(W;\ZZ),A) \cong \cohomol W.
\]
The last isomorphism is an easy consequence of the universal
coefficient theorem. If $\beta$ is the first Betti number of
$W$ we then see that
\[
\MF_W \cong A^{\beta} \times \mbox{Tors},
\]
where $\mbox{Tors}$ is a discrete abelian group of torsion
and that we may therefore identify $\MF_W$ with a compact
abelian Lie group.
Denote the normalized Haar measure on this group by
$\meas W$. Note that $\MF_W$ can also be identified with
$[W,K_A]$, the set of based homotopy classes of maps from
$W$ to the Eilenberg--Mac\,Lane space $K_A = K(A,1)$.

The ``action'' of the theory is defined by a cohomology
class $[\theta]\in H^{n+1}(K_A; U(1))$ by
\[
\begin{array}{ccc}
\MF_W & \longrightarrow & \!\!\!\!\!\!\!\! U(1)\\
\nu & \mapsto & \langle \nu^*([\theta]), [W]\rangle,
\end{array}
\]
 where $[W]$ is the fundamental class of $W$ and
$\langle - , - \rangle$ is the evaluation pairing.
Here we have $\nu^*\colon H^{n+1}(K_A;U(1)) \ra H^{n+1}(W;U(1))$,
thinking of $\nu$ as a homotopy class of maps from $W$ to
$K_{A}$.

If the  action is integrable with respect to the measure
$\meas W$ the DW--invariant of $W$ based
on $[\theta]$ is defined to be
\begin{equation}\label{eq:dw}
Z_A^{[\theta]} (W)= \int_{\nu\in\MF_W}
\langle \nu^*([\theta]),[W]\rangle \dmeas W.
\end{equation}

Our interest in such invariants stems from the following.  In general
a topological quantum field theory (TQFT) is either defined in a
geometrically meaningful way via a non-rigorous path integral or
combinatorially, where the link to the underlying geometry is less
clear. It has been a main goal of the subject for years to bring these
two points of view closer together. Dijkgraaf--Witten TQFT, which
begins with invariants defined by formula (\ref{eq:dw}) where $A$ is
replaced with a (not necessarily abelian) finite group $G$, is
rigourously accessible from both perspectives because the path
integral is a finite sum over $[W,BG] \cong \Hom (\pi_1W,G)$.  With
this in mind one would like to extend Dijkgraaf--Witten theory to
compact Lie groups, but in general the path integral becomes
undefined. For the case of a compact abelian Lie group, however, the
theory can still be approached from both points of view\footnote{For
nonabelian gauge groups we would not get a group structure on $\MF_W$,
hence no Haar measure. The existence of a good measure on $\MF_W$ (and
on other spaces of field configurations to be defined later) is the
main reason that we will look only at abelian groups. Most of the
results that do not involve these measure theoretical problems are
also valid for nonabelian gauge groups.}.

Another reason why Dijkgraaf--Witten theories for continuous groups are
interesting is that they can be viewed as state sum models
in which the set of states over which one must sum to
obtain the invariants is no longer finite or even discrete,
but still finite dimensional.
Because of this, the path integral in these theories is at
an intermediate level of difficulty between the finite
sums of conventional state sum models and the infinite
dimensional integrals that usually occur in
non-topological models.

The original motivation for DW--invariants
\cite{DijkgraafWitten} is that they arise as the partition
functions of Chern--Simons theories with finite gauge group.
The physical states correspond to equivalences classes of
principal $G$--bundles, and Dijkgraaf and Witten show that
$H^3(K(G,1);U(1))$ classifies possible actions. The partition
function is only one aspect of a full TQFT lying behind.
One central feature of TQFTs is locality: a
global invariant can be built up from local contributions.
Locality and the problems associated with patching local
information together make the rigorous construction of
Dijkgraaf--Witten TQFT highly non-trivial. This programme
was carried out by Freed and Quinn in \cite{FreedQuinn}
(and there is related work on $U(1)$--Chern--Simons theory
by Manoliu in \cite{Manoliu}). Turaev \cite{Turaev} has
also recast the pre-path-integral structure arising in
Freed and Quinn's work into an different axiomatic
framework with his homotopy quantum field theories.
This set-up is not specific to dimension 2+1 but works for
any dimension.

It can immediately be seen from (\ref{eq:dw}) that the
DW--invariants only depend on the homotopy
type of the manifold $W$. It is interesting to ask how
good these invariants are as homotopy invariants. From
a purely homotopical point of view locality is rather
unnatural: the homotopy theory of local bits will
overlook aspects of the global homotopy theory. If
one is to proceed to understand the full theory, great
care must be taken to work with respect to prescribed
boundary conditions (on the local pieces) and then
carefully analyse how to fit the pieces together. None
the less the invariants themselves, being homotopy
invariants, should be computable in a more natural way
from the point of view of classical techniques in
algebraic topology.  In this paper we wish to follow
something of a middle road, constructing a theory which
is both simple and natural with regard to homotopy theory,
at the expense of sacrificing full locality for a
restricted version. The restriction is that we will
only allow decompositions along {\em connected}
submanifolds.

We adopt the formalism of Turaev's homotopy quantum field theories
(HQFTs) and begin by recalling some background about these. The
idea of integrating an HQFT to give a TQFT is briefly discussed.
Beginning with a compact abelian Lie group $A$ and an HQFT
of a certain type, we construct a version of
Dijkgraaf--Witten theory based on $A$. We refer to this as {\em
abelian homotopy Dijkgraaf--Witten theory} both to indicate the
link with HQFT and to distinguish our
theory from the Freed--Quinn formulation. Such a thing will consist
of the following assignments.
\begin{itemize}
\item To each closed, oriented $n$--manifold $M$ and $\alpha\in\MF_M$
  we assign a line $L_{M,\alpha}$.
\item To each cobordism $W$ with incoming (resp.\ outgoing) boundary
  $M_0$ (resp.\ $M_1$) and $(\alpha_0,\alpha_1) \in
  \MF_{M_0} \times \MF_{M_1}$ we assign a linear map
\[
\KK W(\alpha_0,\alpha_1) \colon L_{M_0,\alpha_0} \ra
L_{M_1,\alpha_1}.
\]
\end{itemize}
A notable feature is that the construction works in any dimension. We
examine properties of such theories, in particular we prove a
decomposition formula (Theorem~\ref{thm:decomp}) and examine
invariants of products (Theorem~\ref{thm:products}).  We devote the
final section to calculations using both decomposition and product
formulae but also showing how the more familiar combinatorial picture
emerges for explicit calculation. We show for example in
Theorem~\ref{thm:lens} that the DW--invariants with group
$A=U(1)$ separate lens spaces up to homotopy type.


\section{Background on HQFTs}

\subsection{What is an HQFT?}

An HQFT may be seen as an axiomatic formulation of the
``action'' in a TQFT in which the spaces of fields on
a closed $n+1$--manifold $W$ is the set of homotopy
classes of maps from $W$ to some auxiliary space $X$.
Typically $X$ will be an Eilenberg--Mac\,Lane space for
a discrete group and hence the spaces of fields is
related to the moduli space of flat bundles with
connection. This is, in fact, the motivating example
and is a formulation of the ``extended action'' found
in Freed and Quinn's work on Chern--Simons theory for
finite gauge group \cite{FreedQuinn}. HQFTs were defined
by Turaev in \cite{Turaev} (and in a special case in
\cite{BrightwellTurner} and further discussion of the
connection between the two can be found in
\cite{Rodrigues}).

To formulate the theory one considers smooth, oriented,
closed $n$--manifolds and their diffeomorphisms and
cobordisms between such. An $n+1$--dimensional cobordism
(or $n+1$--cobordism for short) is a triple $(M_0,W,M_1)$
where $W$ is a smooth oriented $n+1$--manifold whose
boundary is a disjoint union of $n$--manifolds $M_0$ and
$M_1$ such that $M_1$ has the induced orientation and
$M_0$ the opposite orientation to the induced one. Now
consider all manifolds and cobordisms to come equipped
with characteristic maps, that is to say, maps to some
auxiliary ``background space'' $X$. (Such manifolds and
cobordisms are called $X$--manifolds and $X$--cobordisms
respectively.) Given a $X$--cobordism $(M_0,W,M_1)$
note that by reversing the orientation of $W$ we get a
$X$--cobordism $(M_1,\op{W},M_0)$. It will sometimes be
convenient to write $\op {\sigma}$ for the characteristic
map of this opposite cobordism, where
$\sigma \colon W \ra X$ is the characteristic map of
$(M_0,W,M_1)$.

The key ingredients of an HQFT are assignments as follows.
To each $n$--manifold $M$ with characteristic map
$\gamma\colon M \ra X$ one assigns a finite dimensional
vector space $V_{M,\gamma}$, and to each diffeomorphism one
assigns an isomorphism of these vector spaces. To each
cobordism $(M_0,W,M_1)$ with characteristic map $\sigma$
one assigns a linear map
$V_{M_0,\gamma_0} \ra V_{M_1,\gamma_1}$, where $\gamma_0$
and $\gamma_1$ are the characteristic maps induced on the
boundary. These assignments are subject to a list of
axioms and the reader is asked to consult \cite{Turaev}
for details. Key among the axioms is that the linear maps
associated to $X$--cobordisms are invariant under
homotopies of the characteristic map. It is also worth
noting that Turaev's axiom 1.2.7 has a somewhat special
status. (Here and elsewhere Turaev's axioms refer to the
axioms in \cite[Sect.~1.2]{Turaev}.)
For a general background space it may be
undesirable to impose this axiom as it reduces the theory
to the theory based on an Eilenberg--Mac\,Lane space. In
this paper the background space will be an
Eilenberg--Mac\,Lane space and we will make use of this
axiom.

\subsection{Example}

The following class of examples is due to Turaev, cf.\
\cite[Sect.~1.3]{Turaev}. They are
{\em rank one} in the sense that all vector spaces associated to
$n$--manifolds are one dimensional.

\begin{exe}\label{exe:Turaev}
(Turaev) {\em Primitive cohomological} HQFTs
\end{exe}

Let $X$ be any topological space and let
$\theta\in C^{n+1}(X;U(1))$. For $\gamma \colon M\ra X$
set
\[
L_{M,\gamma} = \CC \{ a\in C_nM \; | \; [a] = [M]\}
/ a \sim\gamma^*\theta(e) b.
\]
In the above $e\in C_{n+1}M$ such that
$\partial e = -a+b$. For $\sigma\colon W\ra X$ define
a homomorphism
\[
E_{W,\sigma}\colon L_{M_0,\gamma_0} \ra  L_{M_1,\gamma_1}
\]
on generators by
\[
a_0 \mapsto \sigma^*\theta (f) a_1,
\]
where $f\in C_{n+1}W$ satisfies
$\partial f = - a_0 + a_1$ and is a representative of
the fundamental class in $H_{n+1}(W,\partial W)$.
Turaev shows that this construction is independent of
any choices and that it indeed gives rise to an HQFT.
Moreover, cohomologous cocycles give equivalent theories.

The following lemma is immediate from the definition just
given.

\begin{lem}\label{lem:inv}
A primitive cohomological HQFT satisfies
\[
E_{\op{W},\op{\sigma}} = E_{W,\sigma}^{-1}.
\]
\end{lem}

In the next section we will restrict to the case where $X$
is an Eilenberg--Mac\,Lane space $K(A,1)$ for a compact
abelian Lie group $A$. Sometimes it will be convenient to
consider group cocycles instead of singular cocycles which
we can do using the fact that the cohomology of the space
$K(A,1)$ is isomorphic to the group cohomology of the
(discrete) group $A$. Recall that a
{\em group $n$--cochain} with coefficients in $U(1)$ is
a function $\omega\colon A^{n} \ra U(1)$ and such
functions form a group $K^n$ under pointwise
multiplication.  The {\em coboundary operator}
$\delta \colon K^n \ra K^{n+1}$ is defined by
\begin{eqnarray*}
\label{groupdelta}
\delta \omega(x_1,\ldots,x_{n+1})&=&
\omega(x_2,\ldots,x_{n+1})
\omega^{-1}(x_1 x_2, x_3, \ldots,x_{n+1})
\omega(x_1,x_2 x_3, x_4, \ldots,x_{n+1}) \\
~&~&\ldots\omega^{(-1)^{n+1}}(x_1,\ldots,x_n x_{n+1})
\omega^{(-1)^{n+2}}(x_1,\ldots,x_n).
\end{eqnarray*}
We have $\delta^2 = 0$ and the group cohomology is
defined as the homology of this cochain complex.
A group $n$--cocycle $\omega$ is {\em normalized} if the
function $\omega$ takes the value $1$ whenever at least
one entry is the identity.

Thus given $\theta \in H^{n+1}(K(A,1);U(1))$ we may choose
a corresponding group cocycle $\omega$. This group cocycle
is not necessarily normalized, but within the cohomology
class of $\theta$ one can always choose a normalized
representative. Conversely, given a group cocycle one
can choose a corresponding singular cocycle representing
the same cohomology class.

\begin{exe}\label{exe:cocyle}
{\em Taking $A=U(1)$ we define 
a 3--cocycle $\theta_k \in C^3(K(U(1),1);U(1))$ and
the corresponding group cocycle
$\omega_k\colon U(1)^3 \ra U(1)$ for any integer $k$.

Noting that $H^1(K(U(1),1); U(1)) = [K(U(1),1),K(U(1),1)]$
pick a $1$--cocycle representing the identity map. Lift this
to a real cochain $\eta\colon K(U(1),1) \ra \RR$. This is
not a cocycle but $\delta \eta$ takes integer values. Now
consider the real valued $3$--cochain
$\eta \cup \delta \eta$ which again is not a cocycle,
however
$\delta (\eta \cup \delta \eta) =
\delta \eta \cup \delta \eta$
which has integer values. We then define
$\theta_k\in C^3(K(U(1),1) ; U(1)) =
\Hom (C_3(K(U(1),1); \ZZ) , U(1))$
by
\[
\theta_k = e^{2\pi i k \eta \cup \delta \eta}.
\]
Note that $\theta_k$ is independent of the lift $\eta$
and is now a cocycle. To define $\omega_k$, let
$g_1, g_2,g_3\in A$ and write
$g_{1} = e^{2\pi i a}, g_{2} = e^{2\pi i b}$ and
$g_{3} = e^{2\pi i c}$ with $0\leq a,b,c < 1$. Then set
\begin{equation}\label{eq:gpcocycle}
\omega_k (g_{1},g_{2}, g_{3}) = e^{2\pi i k a ( b+c - [b+c])},
\end{equation}
where the square bracket means addition modulo $1$.}
\end{exe}

\begin{exe}\label{exe:u1xu1}
{\em When $A=U(1)\times U(1)$, we have the cocycles associated to
the individual $U(1)$ factors, defined above, but we also get a
second type of cocycles. These cocycles are also labelled by an
integer $l \in \ZZ$  and we will call them $\zeta_l$. The
definition of $\zeta_l$ is very similar to that of $\theta_k$.
First, we define $1$--cocycles corresponding to the identity maps
of the first and second factors of $K(U(1)\times U(1),1)\cong
K(U(1),1)\times K(U(1),1)$. Then we lift these to real cochains
$\eta_1,\eta_2$ and we note that the boundaries of these cochains
take integer values. As a consequence, the same is true for the
boundary of the real cochain $\eta_1\cup\delta\eta_2$ and hence we
may define the $3$--cocycles $\zeta_l$ by
\[
\zeta_l = e^{2\pi i l \eta_1 \cup \delta \eta_2}.
\]
To write down the corresponding group cocycles $\psi_l$,
we introduce a similar notation to the one used in formula
(\ref{eq:gpcocycle}). Let $g_1,g_2,g_3\in A=U(1)\times U(1)$
and write $g_{1} = (e^{2\pi i a_1},e^{2\pi i a_2})$,
$g_{2} = (e^{2\pi i b_1}, e^{2\pi i b_2})$ and
$g_{3} = (e^{2\pi i c_1},^{2\pi i c_2})$ with
$0\leq a_i,b_i,c_i < 1$. Then $\psi_l$ is given by
\begin{equation}\label{eq:u1u1gpcocycle}
\psi_l (g_{1},g_{2}, g_{3}) =
e^{2\pi i l a_1 ( b_2+c_2 - [b_2+c_2])},
\end{equation}
where the square bracket means addition modulo $1$ as in
(\ref{eq:gpcocycle}).}
\end{exe}

Clearly, we might also have reversed the roles of the first
and second $U(1)$ factors in example~\ref{exe:u1xu1}.

\subsection{TQFTs and matrix elements}\label{subsec:tqft}

To obtain a TQFT from an HQFT one should perform some kind
of integration. Although this may not be rigourously
defined it is useful to keep it in mind and in the
following brief digression we give an outline of the idea.

Recall that a rank one HQFT assigns a (complex) line
$L_{M,\gamma}$ to each $X$--manifold $(M,\gamma)$. One
should think about the collection of these as a line
bundle $L_M$ over $\Map (M,X)$, the space of fields of
$M$. The Hilbert space associated to $M$ is the space of
sections of this line bundle. The time evolution $U_{W}$
along a cobordism $W$ denoted $U_W$ is defined on a
section $\psi$ of $L_{M_0}$ by
\[
U_W(\psi)(\gamma_1) = \int_{\gamma_0\in\Map(M_0,X)}
K_W(\gamma_0,\gamma_1)(\psi(\gamma_0))d\!\gamma_0,
\]
where $\gamma_1\in\Map(M_1,X)$ and the $
K_W(\gamma_0,\gamma_1)$ are the ``matrix elements'' of
the theory. In this context a matrix element is a linear
map
$K_W(\gamma_0,\gamma_1)\colon
L_{M_0,\gamma_0} \ra L_{M_1,\gamma_1}$
defined by
\[
K_W(\gamma_0, \gamma_1)(x) =
\int_{\sigma\in \Map (W,X; \gamma_0, \gamma_1)}
E_{W,\sigma}(x)d\!\mu,
\]
where $\Map (W,X; \gamma_0, \gamma_1)$ consists of maps
$W\ra X$ agreeing with the given $\gamma_0$ and $\gamma_1$
on the incoming and outgoing boundaries.

The reader should not require much convincing that in
general much of this is ill defined. It is, however, worth
noting that since the homomorphisms $E_{W,\sigma}$ are
homotopy invariant, the ``measure'' $d\!\mu$ needs only
defining on homotopy classes rather than the full
mapping space, which may simplify the situation.

The fundamental property of locality can be expressed in
terms of the matrix elements as follows. Suppose that $W$
can be decomposed along
$M$ as $W=W^\prime \cup_M W^\dprime$. Then locality is
the requirement
\[
K_{W^\prime \cup W^\dprime} (\gamma_0, \gamma_1)(x) =
\int_{\gamma\in \Map (M,X)}
(K_{W^\dprime}(\gamma,\gamma_1) \circ
K_{W^\prime}(\gamma_0,\gamma))(x) d\!\gamma.
\]

Numerical invariants of closed manifolds arise in the
usual way: regard a closed, oriented $n+1$--manifold $W$
as a cobordism from $\emptyset$ to $\emptyset$ in which
case $K_W(\emptyset,\emptyset)(1)\in\CC$ defines
a numerical invariant of $W$.


\section{The definition of abelian homotopy
Dijkgraaf--Witten theory}

We now turn to the central topic of the paper. {\em For the
remainder of the paper $A$ will denote a
compact abelian Lie group and $K_A$ will denote
the Eilenberg--Mac\,Lane space $K(A,1)$.}
The space $K_A$ may be considered as the classifying space
of $A$ regarded as a discrete group.
We note that $A$ is isomorphic to the product of
a torus and a finite abelian group, see e.g.\
\cite[Corollary I.3.7]{BrockerDieck}. We will freely use the
fact that
\[
\cohomol W \cong [W,K_A] \cong \Hom(H_1(W;\ZZ),A),
\]
where the square bracket refers to based homotopy classes
of maps.

Given an $n$--manifold $M$ set
\[
\MF_M = \cohomol M
\]
and similarly for an $n+1$--cobordism $(M_0,W,M_1)$ set
\[
\MF_W = \cohomol W.
\]
These are the ``fields'' of the theory and can be
identified with isomorphism classes of principal
$A$--bundles with flat connection. There is a natural
topology on $\MF_W$ arising from the identification
$ H^1(-;A) \cong \Hom(H_1(-;\ZZ), A)$ which shows that
$\MF_W$ can be identified with the product of a number of
copies of $A$ and a discrete abelian group of torsion.

For any submanifold $M$ of $W$ the inclusion
$i\colon M \ra W$ induces a restiction map
$i^* \colon \MF_W \ra \MF_M$ which we will
denote $r_M^W$.

\begin{lem}\label{lem:rescont}
The restriction map $r_M^W \colon \MF_W \ra \MF_M$ is
continuous.
\end{lem}

\begin{proof}
It suffices to show that composition with the projection
$p$ onto each factor in $\MF_M=A^l \times \mbox{Tors}_M$
is continuous. Let $B$ be such a factor and $Z$ be the
corresponding cyclic group factor in $H_1(M;\ZZ)$ i.e.\
$B = \Hom (Z, A)$. Then
$p\circ r \colon \MF_W =A^k \times \mbox{Tors}_W \ra B$
maps
$(a_1,\ldots,a_k,b_1, \ldots ,b_q)$ to $a_1^{n_1}
\cdots a_k^{n_k}b_1^{m_1} \cdots b_q^{m_q}$,
where the map
\[
Z \ra H_1(M;\ZZ) \ra H_1(W;\ZZ) = \ZZ^k \times \ZZ/p_1
\times \cdots \times \ZZ/p_q
\]
takes $1$ to $(n_1,\ldots , n_k,m_1, \ldots ,m_q)$.
Since multiplication in $A$ is continuous this shows that
$p \circ r$ is continuous.
\end{proof}

For a cobordism $(M_0,W,M_1)$ we will need to consider
fields with prescribed boundary conditions. For a given
pair $(\alpha_0,\alpha_1) \in \MF_{M_0} \times \MF_{M_1}$
of boundary fields we set
\[
\MFr W {\alpha_0} {\alpha_1} =
\{\nu \in \MF_W \mid r_{M_0}^W(\nu) =
\alpha_0 \mbox{ and } r_{M_1}^W(\nu) = \alpha_1 \}.
\]
By Lemma~\ref{lem:rescont}, $\MFr W {\alpha_0} {\alpha_1}$
is a closed, hence compact subset of $\MF_W$ (perhaps
empty).

\subsection{An HQFT--like construction}

Suppose we are given a primitive cohomological HQFT in
dimension $n+1$ with background space $K_A = K(A,1)$.
Let $M$ be an $n$--manifold and let
$\gamma, \gamma^\prime \colon M \ra K_A$.

\begin{prop}
If $\gamma$ is homotopic to $\gamma^\prime$ then
$L_{M,\gamma}$ is canonically isomorphic to
$L_{M,\gamma^\prime}$.
\end{prop}

\begin{proof}
Let $h\colon M\times I \ra K_A$ be a homotopy. Regarding
this as the characteristic map of a cobordism, the HQFT
gives rise to an isomorphism
\[
E_{M\times I,h} \colon L_{M,\gamma} \ra L_{M,\gamma^\prime}.
\]
Given another homotopy $h^\prime\colon M\times I \ra K_A$
consider the map
\[
h\cup \op {h^\prime} \colon M\times I \ra K_A
\]
defined by $h$ on the first half of the cylinder and by
$\op{h^\prime}$ on the second half. This map satisfies
$h\cup \op {h^\prime}|_0 = \gamma $ and
$h\cup \op {h^\prime}_1 = \gamma $ so by the axioms of
HQFTs (and in particular Turaev's axiom 1.2.7 which holds
since our background space is an Eilenberg--Mac\,Lane
space) we have
\[
E_{M\times I , \op{h^\prime}} \circ E_{M\times I , h}
  = E_{M\times I , h \cup \op{h^\prime}} = Id.
\]
Using Lemma~\ref{lem:inv} we conclude
\[
E_{M\times I , h} = E^{-1}_{M\times I , \op{h}}
  = E_{M\times I ,h^\prime}.
\]
Hence the isomorphism given above is independent of the
choice of homotopy finishing the proof.
\end{proof}

This proposition means that given $\alpha\in\MF_M$ we
can define a one dimensional vector space $L_{M,\alpha}$
by identifying the $L_{M,\gamma}$ given by the HQFT using
the canonical isomorphisms above i.e.\ denoting the
isomorphisms above by $\sim$ set
\[
L_{M,\alpha} = \bigoplus_{\{\gamma \mid [\gamma] = \alpha \}}
L_{M,\gamma} / \sim .
\]

Next, given $\alpha_0\in\MF_{M_0}$, $\alpha_1\in\MF_{M_1}$
and $\nu\in\MFr W {\alpha_0}{\alpha_1}$ we wish to define
\[
E_{W,\nu}\colon L_{M_0,\alpha_0} \ra L_{M_1,\alpha_1}.
\]
Suppose that
$\sigma\colon W\ra K_A, \gamma_0\colon M_0\ra K_A$ and
$\gamma_1\colon M_1 \ra K_A$ are maps representing $\nu$,
$\alpha_0$ and $\alpha_1$ respectively. Suppose moreover
that
$\gamma_0 = \sigma |_{M_0}$ and $\gamma_1 = \sigma |_{M_1}$.
Courtesy of the HQFT we have a map
\[
E_{W,\sigma}\colon L_{M_0,\gamma_0} \ra L_{M_1,\gamma_1},
\]
which induces a map
\[
L_{M_0,\alpha_0} \ra L_{M_1,\alpha_1}.
\]

\begin{prop}\label{prop:homotopyindependence}
The induced map above depends only on the homotopy class of $\sigma$.
\end{prop}

\begin{proof}
Let $\sigp$ be another choice with
$\sigp |_{M_0} = \gpzero$ and $\sigp |_{M_1} = \gpone$.
In order to prove the proposition we must
show that the following diagram commutes.
\[
\xymatrix{
L_{M_0,\gamma_0} \ar[r]^{E_{W,\sigma}}
\ar[d]_{c_{\gamma_0,\gpzero}}
 & L_{M_1,\gamma_1} \ar[d]^{c_{\gamma_1,\gpone}} \\
L_{M_0,\gpzero} \ar[r]_{E_{W,\sigp}} & L_{M_1,\gpone}
}
\]
(the vertical maps are the canonical isomorphisms above).
Let $H$ be a homotopy from $\sigma$ to $\sigp$ and let
$h_0=H|_{M_0\times I}$ and $h_1 = H|_{M_1\times I}$.
Note that $h_0$ is a homotopy from $\gamma_0$ to
$\gpzero$ and that $h_1$ is a homotopy from $\gamma_1$
to $\gpone$. Consider
\[
W^\prime = (M_0 \times I) \cup_{M_0} W \cup_{M_1}
(M_1\times I)
\]
and let $g\colon W^\prime \ra K_A$ be defined by
$g = h_0 \cup \sigp \cup h_1^{-1}$. Using the HQFT and its
properties we get
\[
E_{W^\prime, g} = E_{M_1\times I,h_1^{-1}} \circ
  E_{W,\sigp} \circ E_{M_0\times I, h_0}
= c_{\gamma_1,\gpone}^{-1} \circ  E_{W,\sigp}
  \circ c_{\gamma_0,\gpzero}.
\]
Hence in order to show that the diagram above commutes we
need to show that $E_{W,\sigma} = E_{W^\prime,g}$.

Define $H_0\colon M_0\times I \times I \ra K_A$ by
$H_0(x,s,t) = h_0(x,s(1-t))$ and define
$H_1\colon M_1\times I \times I \ra K_A$ by
$H_1(x,s,t) = h_1(x,(1-s)(1-t))$. Note that $H_0$ provides
a homotopy between $h_0$ and $\kappa_{\gamma_0}$, where
the latter is defined by
$\kappa_{\gamma_0}(x,s) = \gamma_0 (x)$. Similarly $H_1$
provides a homotopy between $h_1^{-1}$ and
$\kappa_{\gamma_1}$. Now define a map
\[
\MH \colon W^\prime \times I = ((M_0 \times I)
 \cup_{M_0} W \cup_{M_1} (M_1\times I))\times I \ra K_A
\]
by $\MH = H_0 \cup H^{-1} \cup H_1$. This map provides a
homotopy between $g$ and
$f= \kappa_{\gamma_0} \cup \sigma \cup \kappa_{\gamma_1}$.
Moreover, it is readily checked that on the boundary $\MH$
is $\kappa_{\gamma_0} \sqcup \kappa_{\gamma_1}$ which is
independent of $t$. Thus $\MH$ provides a homotopy rel
$\partial W^\prime$ from $g$ to $f$ and hence by the
properties of an HQFT (see Turaev's axiom 1.2.8) we have
\[
E_{W^\prime,g} = E_{W^\prime, f}.
\]
There is a diffeomorphism $T\colon W^\prime \ra W$ making
the following diagram commute.
\[
\xymatrix
{
W^\prime \ar[rr]^T \ar[dr]_f & & W \ar[dl]^\sigma \\
& K_A &
}
\]
Hence, by Turaev's axiom 1.2.4 the following diagram
commutes.
\[
\xymatrix
{
L_{M_0,\gamma_0} \ar[r]^{E_{W^\prime,f}} \ar[d]_{id} &
L_{M_1,\gamma_1} \ar[d]^{id}\\
L_{M_0,\gamma_0} \ar[r]_{E_{W,\sigma}} & L_{M_1,\gamma_1}
}
\]
Thus $ E_{W,\sigma} = E_{W^\prime, f} = E_{W^\prime,g}$
which finishes the proof.
\end{proof}

Given $\alpha_0\in\MF_{M_0}$, $\alpha_1\in\MF_{M_1}$ and
$\nu\in\MFr W {\alpha_0}{\alpha_1}$ the above proposition
shows that we have a well-defined map
\[
E_{W,\nu}\colon L_{M_0,\alpha_0} \ra L_{M_1,\alpha_1}
\]
defined by $E_{W,\nu} = E_{W,\sigma}$, where $\sigma$ is
any representative of the class $\nu$. As a corollary of
Lemma~\ref{lem:inv} we have
\begin{equation}\label{eq:inv}
E_{\op{W}, \nu} = E^{-1}_{W,\nu}.
\end{equation}
If $W$ is a closed manifold and $\sigma\in \MF_W$ then
$\sigma\colon W \ra K_A$ may be considered as the
classifying map of a principal $A$--bundle. The invariant
$E_{W,\sigma}(1)\in\CC^\times$ should correspond to
Turaev's invariant $\tau_{\MC} (W,\sigma)$
constructed (via surgery) in \cite{turaevhqft3},
where $\MC$ is the modular
$A$--category constructed from the cocycle $\theta$.

\subsection{Measures and abelian homotopy
Dijkgraaf--Witten theory}\label{subsec:measure}

In order to construct our analogue of matrix elements
we need to integrate and in order to integrate we need
to put measures on our spaces of fields.

We have identified $\MF_W$ as the product of a number of
copies of $A$ and a discrete abelian group of torsion.
Thus we can equip $\MF_W$ with the normalized Haar measure
which we denote by $\meas W$. Since $A$ is the product of
a torus and a finite abelian group we have that $\MF_W$ is
also the product of a torus $T$ and a finite abelian group
$B$ and it follows by the defining properties of the Haar
measure that the Haar measure on $\MF_W$ is nothing but
the product of the Lebesgue measure on $T$ and the
counting measure on $B$ (normalized). For details on the
Haar measure and the associated Haar integral we refer to
\cite[Chap.~VIII]{Knapp}, \cite{Nachbin} and
\cite[Chap.~6]{Pedersen}. Let us just remark here that
the left invariant Haar measures on a Lie group $G$
(which all differ by a scalar) are Borel measures, i.e.\
they are measures on the $\sigma$--algebra of all Borel
sets of $G$. Moreover, if $G$ is abelian or compact
then the normalized left and right invariant Haar measures
coincide and are just called the normalized Haar measure
on $G$, see e.g.\ \cite[Corollary 8.31]{Knapp} or
\cite[p.~81]{Nachbin}.

Let us now define a measure on each of the spaces
$\MFr W {\alpha_0}{\alpha_1}$. It is tempting to define
this measure using the restriction of the measure on
$\MF_W$, but we will not do this because it would yield
the zero measure whenever $\MFr W {\alpha_0}{\alpha_1}$
has measure zero in $\MF_W$. Instead, we will use the
group structure of $\MF_W$ to define a normalized
measure on each of the $\MFr W {\alpha_0}{\alpha_1}$.
Firstly, note that $\MFr W 0 0$ is a subgroup of $\MF_W$
and being closed it is in fact a compact Lie subgroup of
$\MF_W$, hence we can endow $\MFr W 0 0 $ with its
normalized Haar measure, which we will denote by
$\measr W 0 0 $.

The set $\MFr W {\alpha_0}{\alpha_1}$ is either empty or
a coset of $\MFr W 0 0$, hence the measure $\measr W 0 0 $
induces a normalized measure $\measalpha$ on
$\MFr W {\alpha_0}{\alpha_1}$, namely $\measalpha$ is
nothing but the image measure of $\measr W 0 0 $ under
the translation with an arbitrary element of
$\MFr W {\alpha_0}{\alpha_1}$. That is
\[
\measalpha (S) = \measr W 0 0  (S-\nu)
\]
for any $\nu \in \MFr W {\alpha_0}{\alpha_1}$. (By translation
invariance of the Haar measure this does not depend on the choice
of $\nu$.) We use here and in the following the standard
definition of image measures. That is, given a measurable function
$f \colon X \ra Y$ between two measurable spaces and given a
measure $\mu$ on $X$, we define the image measure of $\mu$ under
$f$ to be the measure $\nu$ on $Y$ given by $\nu(S) =
\mu(f^{-1}(S))$ for any measurable subset $S$ of $Y$. It is then a
standard result in integration theory that if $g \colon Y \ra \CC$
is measurable, then
\begin{equation}\label{eq:intimage}
\int_{y \in Y} g(y) d\!\nu = \int_{x \in X} g(f(x)) d\!\mu
\end{equation}
in the sense that if one of the integrals exists, then so
does the other and the two integrals are equal.

If the measure of $\MFr W 0 0$ in $\MF_W$ is non-zero,
then it follows from the defining properties of the
Haar measure that the measure $\measr W 0 0$ is nothing
but the normalization of the measure obtained by
restriction, hence the same is true for the measures
$\measalpha$ in that case.

We now define homomorphisms
\[
\KK W(\alpha_0,\alpha_1) \colon L_{M_0,\alpha_0}
  \ra L_{M_1,\alpha_1}
\]
by
\begin{equation}\label{eq:matrix}
\KK W(\alpha_0,\alpha_1)(a_0) =
\int_{\nu\in\MFr W {\alpha_0}{\alpha_1}}
E_{W,\nu}(a_0) \dmeasalpha .
\end{equation}
By convention, if
$\MFr W {\alpha_0}{\alpha_1} = \emptyset$, we take
$\KK W(\alpha_0,\alpha_1)$ to be the zero map.
Of course, for the integral in the definition to make
sense, we must insist that the function
$\MFr W {\alpha_0} {\alpha_1} \ra L_{M_1,\alpha_1}$
given by $\nu\mapsto E_{W,\nu}(a_0)$ is integrable. If
all such functions are indeed integrable (and hence the
$\KK W (\alpha_0,\alpha_1)$ are defined) we will refer to
the defining primitive cohomological HQFT as
{\em integrable}.

So, we start with a primitive cohomological HQFT based on
a cocycle $\theta\in C^{n+1}(K_A;U(1))$ and define the 
associated {\em abelian homotopy Dijkgraaf--Witten theory}
to consist of the assignments above. Namely,
\begin{itemize}
\item to each closed, oriented $n$--manifold $M$ and 
  $\alpha\in\MF_M$ we assign the line $L_{M,\alpha}$,
\item to each $n+1$--cobordism $(M_0,W,M_1)$ and 
  $(\alpha_0,\alpha_1) \in \MF_{M_0} \times \MF_{M_1}$ 
  we assign the linear map
\[
\KK W(\alpha_0,\alpha_1) \colon L_{M_0,\alpha_0} \ra
L_{M_1,\alpha_1}.
\]
\end{itemize}

\subsection{Invariants of closed manifolds}

A closed oriented $n+1$--manifold may be regarded as a
cobordism from $\emptyset$ to $\emptyset$ and thus
$\KK W (\emptyset,\emptyset)$ is a map from $\CC$ to $\CC$.
The DW--invariant of $W$ is defined to be
the image of $1$ i.e.\
\[
K_W(\emptyset,\emptyset)(1) = \int_{\nu\in\MF_W}
E_{W,\nu}(1) \dmeas W.
\]
Note that for a given $\sigma\colon W \ra K_A$ the map
$E_{W,\sigma}\colon \CC \ra \CC$ is given by
$1 \mapsto \sigma^*\theta (f)$, where $f\in C_{n+1}W$ is
a fundamental cycle for $W$. It follows that
$E_{W,\sigma} (1)$ is a function of the cohomology class
of $\theta$ only. Thus writing $\langle - , - \rangle$
for the evaluation map
$H^{n+1}(W;U(1)) \otimes H_{n+1}(M;\ZZ) \ra U(1)$ we get
\[
K_W(\emptyset,\emptyset)(1) = \int_{\nu\in\MF_W}
\langle \nu^*([\theta]), [W] \rangle \dmeas W.
\]
Writing $Z_A^{[\theta]} (W)$ for
$K_W(\emptyset,\emptyset)(1)$
this is the expression (\ref{eq:dw}) in the introduction.
If we are given a group cocycle we will sometimes use
the notation $Z_A^{[\omega]}$ instead.

Note that if $W$ and $W^\prime$ are two closed, oriented
$n+1$--manifolds then
\[
Z_A^{[\theta]} (W \sqcup W^\prime) = Z_A^{[\theta]}
(W)Z_A^{[\theta]} (W^\prime).
\]

\begin{exe}\label{ex:spheres}
Spheres.
{\em Let $\theta \in C^{n+1}(K_A;U(1))$ be a cocycle.
For $n>0$ we have
$\MF_{S^{n+1}} = \cohomol {S^{n+1}} = \{ 0 \}$
and so
\[
 Z_A^{[\theta]} (S^{n+1}) = \int_{\nu\in\MF_{S^{n+1}}}
\langle \nu^*([\theta]), [S^{n+1}] \rangle
\dmeas {S^{n+1}} = \langle 0^*([\theta]), [S^{n+1}]
\rangle =  \langle 1, [S^{n+1}] \rangle = 1.
\]
For $n=0$ suppose we have a corresponding group cocycle
$\omega \colon A \ra U(1)$. Noting that $\MF_{S^1} = A$
we have
\[
 Z_A^{[\theta]} (S^{1}) = \int_{\nu\in\MF_{S^{1}}}
\langle \nu^*([\theta]), [S^{1}] \rangle \dmeas {S^{1}}
 = \int_{\nu\in\MF_{S^{1}}}
   \langle \theta , \nu_*[S^{1}] \rangle \dmeas {S^{1}}
 = \int_{a\in A } \omega (a) \dmeas {S^{1}}.
\]
Note here that the $1$--cocycle $\omega$ is just a
$1$--dimensional representation of $A$. Hence the integral
equals $0$ unless $\omega$ is the trivial representation,
in which case the integral equals $1$.}
\end{exe}


\section{Properties of abelian homotopy
Dijkgraaf--Witten theory}

\subsection{Decompositions}

In this section we discuss the restricted version of
locality satisfied by abelian homotopy Dijkgraaf--Witten
theory. Suppose that we can decompose an
$n+1$--cobordism $(M_0,W,M_1)$ into two pieces $W^\prime$
and $W^\dprime$ along a connected $n$--manifold $M$.
Given such a decomposition and given a pair
$(\alpha_0,\alpha_1) \in\MF_{M_0} \times \MF_{M_1}$ we
define the space of {\em supporting fields} to be
\[
\suppfields M {\alpha_0}{\alpha_1} =
\{ \alpha \in \MF_M \;|\; \MFr {W^\prime} {\alpha_0}{\alpha}
 \times \MFr {W^\dprime} {\alpha}{\alpha_1}
   \neq \emptyset \}.
\]
Note that this depends on the decomposition. In this subsection we
construct a measure $\suppmeas M {\alpha_0}{\alpha_1}$ on the space of
supporting fields and we prove the following theorem.

\begin{thm}\label{thm:decomp}
Suppose we can decompose $W$ as
$W=W^\prime \cup_M W^\dprime$, where $M$ is a connected
$n$--manifold and $W^\prime \cap W^\dprime = M$. Then
for $\alpha_0\in\MF_{M_0}$ and $\alpha_1\in\MF_{M_1}$ we
have
\[
\KK W (\alpha_0,\alpha_1) (x) = \int_{\alpha\in\suppfields M
  {\alpha_0}{\alpha_1}} \KK {W^\dprime} (\alpha, \alpha_1) \circ \KK
  {W^\prime} (\alpha_0,\alpha)(x) \dsuppmeas M
  {\alpha_0}{\alpha_1}.
\]
\end{thm}

Before proving this theorem we need to construct the measure
$\suppmeas M {\alpha_0}{\alpha_1}$ on $\suppfields M
{\alpha_0}{\alpha_1}$. The connectedness of $M$ which will be
essential in the proof of this theorem will not be needed for the
construction of the measure so to begin with we will not
assume that $M$ is connected.

By Lemma~\ref{lem:rescont} we have a continuous (restriction) map
$r_M^W\colon \MF_W \ra \MF_M$, which restricts to a continuous
surjection $r=r^{\alpha_{0},\alpha_{1}} \colon \MFr W
{\alpha_0}{\alpha_1} \ra \suppfields M {\alpha_0}{\alpha_1} $.  To
see that $r$ is surjective, apply the Mayer--Vietoris sequence for
the triad $(W;W^\prime,W^\dprime)$, i.e.\ the exact sequence
\begin{equation}\label{eq:MVS}
\cdots \lra \tilde{H}^0 (M) \lra H^1(W) \stackrel{a}{\lra}
H^1(W^{\prime}) \oplus H^1(W^{\dprime}) \stackrel{b}{\lra} 
H^1(M) \lra \cdots,
\end{equation}
where
$a(\nu)=(r_{W'}^W(\nu),r_{W''}^W(\nu))$ and
$b(\nu^\prime, \nu^\dprime) = r_M^{W'} (\nu^\prime) -
r_M^{W^\dprime} (\nu^\dprime)$
(all cohomology groups having coefficients in $A$).
In particular, $\suppfields M
{\alpha_0}{\alpha_1}$ is a closed subset of $\MF_M$. Given $\alpha
\in \suppfields M {\alpha_0}{\alpha_1}$ we will denote
$r^{-1}(\alpha) \subset \MFr W {\alpha_0}{\alpha_1}$ by $\MFrr W
{\alpha_0}\alpha {\alpha_1}$.  We note that $\MFrr W 0 0 0$ is a
compact Lie subgroup of $\MFr{W}{0}{0}$.

As with the measures on the spaces
$\MFr W {\alpha_0}{\alpha_1}$ it turns out that one should
not take the measure induced from the obvious inclusion
(in this case into $\MF_{M}$). To begin let us instead
note that $\suppfields M {0}{0}$ is a Lie subgroup of
$\MF_M$ (as a closed subgroup of $\MF_M$). We therefore
let $\suppmeas M {0}{0}$ be the normalized Haar measure on
$\suppfields M {0}{0}$. Next observe that
$\suppfields M {\alpha_0}{\alpha_1}$ is either empty or a
coset of $\suppfields M {0}{0}$, hence we can (similarly to
the construction of $\measalpha$ in
Sect.~\ref{subsec:measure}) define
$\suppmeas M {\alpha_0}{\alpha_1}$ to be the image measure
of $\suppfields M {0}{0}$ under the translation by any
element of $\suppfields M {\alpha_0}{\alpha_1}$.

Thinking about this slightly differently let 
$\pi=r^{0,0} \colon  \MFr W {0}{0} \ra 
\suppfields M {0}{0} $, which is a surjective Lie group
homomorphism, inducing a Lie group isomorphism
$\bar{\pi} \colon \MFr W {0}{0}/\MFrr{W}{0}{0}{0} \ra
\suppfields M {0}{0} $. Let $\bar{\mu}$ be the normalized
Haar measure on the quotient
$\MFr W {0}{0}/\MFrr{W}{0}{0}{0}$. Then $\bar{\mu}$ is
also the image measure of $\measr W 0 0 $ under the
canonical projection and $\suppmeas M {0}{0}$ is the
image measure of $\bar{\mu}$ under $\bar{\pi}$ and also
the image measure of $\measr W 0 0 $ under $\pi$.
We use here the obvious fact that if $f\colon G \ra H$ is
a surjective Lie group homomorphism and if $\mu_G$ and
$\mu_H$ are the normalized left invariant Haar measures
on respectively $G$ and $H$, then $\mu_H$ equals the
image measure of $\mu_G$ under $f$.

By (\ref{eq:intimage}) we then get
\begin{equation}\label{eq:intid1}
\int_{\suppfields M {\alpha_0}{\alpha_1}} f
d\!\suppmeas M {\alpha_0}{\alpha_1}
= \int_{\nu \in \MFr W {0}{0}}
f(\pi(\nu)+\rho_{\alpha_0,\alpha_1}) d\!\mu_{W}^{0,0},
\end{equation}
for an integrable function $f$ on
$\suppfields M {\alpha_0}{\alpha_1}$,
where $\rho_{\alpha_0,\alpha_1}$ is an arbitrary element
of $\suppfields M {\alpha_0}{\alpha_1}$. Moreover, if $f$
is an integrable function on $\MFr W {0}{0}$ we have
\begin{equation}\label{eq:intid2}
\int_{\MFr W {0}{0}} f d\!\measr W 0 0
=\int_{p(g) \in \MFr W {0}{0}/\MFrr{W}{0}{0}{0}}
 \left( \int_{h \in \MFrr{W}{0}{0}{0}}
 f(g+h) d\!\mu_{W}^{0,0,0} \right) d\!\bar{\mu},
\end{equation}
where $p$ is the canonical projection and
$\mu_{W}^{0,0,0}$ is the normalized Haar measure on
$\MFrr{W}{0}{0}{0}$. We write things additively since
we deal with abelian groups. The identity
(\ref{eq:intid2}) simply follows by noting that both
sides define a normalized integral which is
left-invariant on the class of continuous functions (see
also \cite[Proposition I.5.16]{BrockerDieck} and
\cite[Theorem 8.36]{Knapp} for a more general
result).

Before proving the decomposition theorem we need one more
result. To establish this result we must assume that $M$
is connected. Let $a \colon \MF_W \ra \MF_{W'} \times \MF_{W''}$ be
the continuous restriction map from the Mayer--Vietoris sequence
(\ref{eq:MVS}).

\begin{lem}\label{lem:ident}
Assume that $M$ is connected. The map $a$ restricts to a
bijection
\[
\MFrr{W}{\alpha_0}{\alpha}{\alpha_1} \cong
\MFr{W'}{\alpha_0}{\alpha} \times
\MFr{W^\dprime}{\alpha}{\alpha_1}
\]
for any $\alpha_0\in\MF_{M_0}$, $\alpha_1\in\MF_{M_1}$
and $\alpha \in \suppfields M {\alpha_0}{\alpha_1}$.
In particular, we have a Lie group isomorphism
\[
\MFrr{W}{0}{0}{0} \cong
\MFr{W'}{0}{0} \times \MFr{W^\dprime}{0}{0}.
\]
\end{lem}

\begin{proof}
Let $\alpha_0\in\MF_{M_0}$ and $\alpha_1\in\MF_{M_1}$ be
fixed. The map $a$ clearly maps
$\MFrr{W}{\alpha_0}{\alpha}{\alpha_1}$ into
$\MFr{W'}{\alpha_0}{\alpha} \times
\MFr{W^\dprime}{\alpha}{\alpha_1}$.
Consider the Mayer--Vietoris sequence (\ref{eq:MVS}).
Since $M$ is connected,  $a$ injects. Assume
$(\nu',\nu'') \in \MFr{W'}{\alpha_0}{\alpha}
\times \MFr{W^\dprime}{\alpha}{\alpha_1}$.
Then $r_M^{W^\prime}(\nu')=\alpha=r_M^{W^\dprime}(\nu'')$
so
$(\nu',\nu'') \in \Ker (b)$. Hence there exists a
$\nu \in \MF_W$ such that $a(\nu)=(\nu',\nu'')$, and
by the very definition of
$\MFrr{W}{\alpha_0}{\alpha}{\alpha_1}$ we see that
$\nu \in \MFrr{W}{\alpha_0}{\alpha}{\alpha_1}$.
\end{proof}

The bijections in the above lemma will all be denoted by
$a$. We now prove Theorem~\ref{thm:decomp}.

\begin{proof} (of Theorem~\ref{thm:decomp})\\
There is only something to prove in case
$\MFr{W}{\alpha_0}{\alpha_1}$ is nonempty, so this we
assume in what follows. Let us start by introducing
some notation. The statement of the theorem is for the
linear maps $\KK W (\alpha_{0},\alpha_{1})$,
but we can  choose fixed basis vectors in the
lines associated with $M_0,M_1$ and $M$ and then
replace these linear maps by their matrix elements
$k_W, k_{W'}, k_{W^\dprime}$, which are just functions
of the boundary configurations $\alpha_0$, $\alpha$ and
$\alpha_1$. Similarly, we introduce the notation
$e_W$, $e_{W'}$ and $e_{W^\dprime}$ for the matrix
elements of the linear maps $E_{W,\sigma}$ that are
integrated to give the maps $\KK W$. This means $e_W$ is
a function on $\MF_W$ and analogously for $e_{W'}$ and
$e_{W^\dprime}$. With this notation, we have
(letting $\rho_{\alpha_0,\alpha_1}$ be an arbitrary element of
$\MFr W {\alpha_0}{\alpha_1}$)
\begin{eqnarray*}
k_W(\alpha_0,\alpha_1) &=&
\int_{\nu \in \MFr W {\alpha_0}{\alpha_1}} e_W (\nu)
d\!\mu_{W}^{\alpha_0,\alpha_1} =
\int_{\nu \in \MFr W {0}{0}}
e_W(\nu+\rho_{\alpha_0,\alpha_1}) d\!\mu_{W}^{0,0} \\
&=& \int_{p(\nu) \in \MFr{W}{0}{0}/\MFrr{W}{0}{0}{0}}
\left(\int_{\sigma \in \MFrr{W}{0}{0}{0}}
e_W (\sigma+\nu+\rho_{\alpha_0,\alpha_1})
d\!\mu_{W}^{0,0,0} \right) d\!\bar{\mu},
\end{eqnarray*}
where the final equality follows by (\ref{eq:intid2}).
Next we apply our Lie group isomorphism $a$ from
Lemma~\ref{lem:ident} to get
\begin{eqnarray*}
&&\int_{\sigma \in \MFrr{W}{0}{0}{0}}
e_W (\sigma+\nu+\rho_{\alpha_0,\alpha_1})
d\!\mu_{W}^{0,0,0} \\
 &&\hspace{.5in} = \int_{(\sigma',\sigma'') \in
\MFr {W^\prime} {0}{0} \times \MFr {W^\dprime} {0}{0}}
e_W (a^{-1}(\sigma',\sigma'')+\nu+\rho_{\alpha_0,\alpha_1})
d\!\mu_{W'}^{0,0} \oplus \mu_{W^\dprime}^{0,0}
\end{eqnarray*}
noting that the product measure
$\mu_{W'}^{0,0} \oplus \mu_{W^\dprime}^{0,0}$ is the
normalized Haar measure on the product Lie group
$\MFr {W^\prime} {0}{0} \times \MFr {W^\dprime} {0}{0}$,
hence the image measure of $\mu_W^{0,0,0}$ under $a$.
Using the map
$a \colon \MF_W \ra \MF_{W'} \times \MF_{W^\dprime}$
to write
$(\nu',\nu'')=a(\nu) \in \MFr {W^\prime} {0}{\alpha_{\nu}}
\times \MFr {W^\dprime} {\alpha_{\nu}}{0}$ for
$\nu \in \MF_W^{0,0}$ and
$(\rho_{\alpha_0,\alpha_1}^\prime,
\rho_{\alpha_0,\alpha_1}^\dprime)
=a(\rho_{\alpha_0,\alpha_1}) \in
\MFr {W^\prime} {\alpha_0}{\beta} \times
\MFr {W^\dprime} {\beta}{\alpha_1}$,
where $\beta = r_M^W(\rho_{\alpha_0,\alpha_1})$
and $\alpha_{\nu}=r_M^W(\nu)$, we get
\[
e_W (a^{-1}(\sigma',\sigma'')+\nu+\rho_{\alpha_0,\alpha_1})
= e_{W'}(\sigma'+\nu'+\rho_{\alpha_0,\alpha_1}')
  e_{W^{\dprime}}(\sigma^{\dprime}+\nu^{\dprime}
     + \rho_{\alpha_0,\alpha_1}^{\dprime})
\]
by the HQFT gluing property, hence
\begin{eqnarray*}
k_W(\alpha_0,\alpha_1) &=& \int_{p(\nu) \in
\MFr{W}{0}{0}/\MFrr{W}{0}{0}{0}} \left(
\int_{\sigma' \in \MFr{W'}{0}{0}}
e_{W'}(\sigma'+\nu'+\rho_{\alpha_0,\alpha_1}')
d\!\mu_{W'}^{0,0}\right) \\
 && \hspace{.5in} \times \left(
\int_{\sigma^{\dprime} \in \MFr{W^{\dprime}}{0}{0}}
e_{W^{\dprime}}(\sigma^{\dprime}+\nu^{\dprime}
 + \rho_{\alpha_0,\alpha_1}^{\dprime})
d\!\mu_{W^{\dprime}}^{0,0} \right) d\!\bar{\mu},
\end{eqnarray*}
by Fubini's theorem. Here
\[
\int_{\sigma' \in \MFr{W'}{0}{0}}
e_{W'}(\sigma'+\nu'+\rho_{\alpha_0,\alpha_1}')
d\!\mu_{W'}^{0,0} =
\int_{x \in \MFr{W'}{\alpha_0}{\alpha_{\nu}+\beta}}
e_{W'}(x) d\!\mu_{W'}^{\alpha_0,\alpha_{\nu}+\beta}
= k_{W'}(\alpha_0,\alpha_{\nu} + \beta)
\]
and
\[
\int_{\sigma^\dprime \in \MFr{W^\dprime}{0}{0}}
e_{W^\dprime}(\sigma^{\dprime} +\nu^{\dprime} +
\rho_{\alpha_0,\alpha_1}^{\dprime})
d\!\mu_{W^\dprime}^{0,0} =
\int_{x \in \MFr{W^\dprime}{\alpha_{\nu}+\beta}{\alpha_1}}
e_{W^\dprime}(x)
d\!\mu_{W^\dprime}^{\alpha_{\nu}+\beta,\alpha_1}
= k_{W^\dprime}(\alpha_{\nu}+\beta,\alpha_1).
\]
Therefore, since $\alpha_{\nu}=\bar{\pi}(p(\nu))$,
\[
k_W(\alpha_0,\alpha_1)=
\int_{p(\nu) \in \MFr{W}{0}{0}/\MFrr{W}{0}{0}{0}}
k_{W'}(\alpha_0,\beta+\bar{\pi}(p(\nu)))
k_{W^{\dprime}}(\beta+\bar{\pi}(p(\nu)),\alpha_1)
d\!\bar{\mu}.
\]
By (\ref{eq:intimage}) and the remarks above
(\ref{eq:intid1}) we then get
\begin{eqnarray*}
k_W(\alpha_0,\alpha_1) &=&
\int_{\alpha \in \suppfields{M}{0}{0}}
k_{W'}(\alpha_0,\beta+\alpha)
k_{W^{\dprime}}(\beta+\alpha,\alpha_1)
d\!\suppmeas{M}{0}{0} \\
 &=& \int_{\alpha \in \suppfields{M}{\alpha_0}{\alpha_1}}
k_{W'}(\alpha_0,\alpha)k_{W^{\dprime}}(\alpha,\alpha_1)
d\!\suppmeas{M}{\alpha_0}{\alpha_1}
\end{eqnarray*}
which is the desired result.
\end{proof}

We end this section with an important corollary to
Theorem~\ref{thm:decomp}.

\begin{cor}
In the set up of \,{\em Theorem~\ref{thm:decomp}}
suppose that $n>1$ and that $H_1(M;\ZZ) =\{0\}$. Then
either $\KK W (\alpha_0,\alpha_1)$ is trivial or
\[
\KK W (\alpha_0, \alpha_1) =
\KK {W^\dprime} (0, \alpha_1) \circ
\KK {W^\prime} (\alpha_0, 0).
\]
\end{cor}

\begin{proof}
Follows immediately from the fact that
$\suppfields M {\alpha_0}{\alpha_1} \subset \MF_M = \{0\}$.
\end{proof}

\subsection{Connected sums}

The decomposition theorem of the previous section allows
us to calculate invariants of connected sums. First we
need the following.

\begin{lem}\label{lem:disks}
If $D$ is an $n+1$--disk with ingoing boundary sphere
then
\[
\KK {\op{D}} (\emptyset , 0 ) \circ \KK D (0,\emptyset)
= \Id.
\]
\end{lem}

\begin{proof}
Let $a$ be a (representative of a) generator of
$L_{S^n,0}$ and note that since $D$ is contractible
$\MFr D 0 \emptyset = \{ 0 \}$. Thus
$\KK D ( 0 , \emptyset) (a) = E_{D,\sigma}(a)$,
where $[\sigma] = 0$ and similarly
$\KK {\op{D}} (\emptyset , 0 )(1) = E_{\op{D},\sigma}(1)$.
Thus using Lemma~\ref{lem:inv} we have
\[
\KK {\op{D}} (\emptyset , 0 ) (\KK D ( 0 , \emptyset) (a))
= \KK {\op{D}} (\emptyset , 0 )(E_{D,\sigma}(a)) =
  E_{\op{D},\sigma}(E_{D,\sigma}(a)) = a.
\]
\end{proof}

Now for the result on connected sums.

\begin{prop}\label{prop:connectedsum}
If $W^\prime$ and $W^\dprime$ are closed, oriented
connected $n+1$--manifolds then
\[
Z^{[\theta]}_A (W^\prime \# W^\dprime )
= Z^{[\theta]}_A (W^\prime
) Z^{[\theta]}_A (W^\dprime ).
\]
\end{prop}

\begin{proof}
Let $(\emptyset, V^\prime, S^n)$ be the $n+1$--cobordism
obtained from $W^\prime$ by removing an $n+1$--disk $D$
(creating a new outgoing boundary component), and similarly
let $(S^n,V^\dprime, \emptyset)$ be the cobordism obtained
from $W^\dprime$ again by removing an $n+1$ disk (this
time creating a new incoming boundary component). We can
then write
$W^\prime \# W^\dprime = V^\prime \cup_{S^n} V^\dprime$.

Also note that $W^\prime = V^\prime \cup_{S^n} D$. For this
decomposition observe that
$\suppfields {S^n} \emptyset \emptyset = \{ 0 \}$.
This is immediate for $n>1$ and for $n=1$ we note that
$V^\prime$ has the homotopy type of a wedge of circles and
that the restriction map to $\MF_{S^1}$ is given by a
commutator map which is trivial since $A$ is abelian and
thus $\MFr {V^\prime} \emptyset \alpha$ is only non-empty
when $\alpha = 0$. Using this we see
\[
Z^{[\theta]}_A(W^\prime) =
\KK {V^\prime \cup_{S^n} D} (\emptyset,\emptyset) (1) =
\KK D ( 0 , \emptyset ) \circ
\KK {V^\prime}(\emptyset , 0 ) (1) .
\]
Similarly
$Z^{[\theta]}_A(W^\dprime)=
\KK {V^\dprime} ( 0 , \emptyset ) \circ
\KK {\op{D}}(\emptyset , 0 ) (1) $.
Thus by applying Theorem~\ref{thm:decomp} we have
\begin{eqnarray*}
Z^{[\theta]}_A(W^\prime \# W^\dprime) &=&
\KK {V^\prime \cup_{S^n} V^\dprime} (\emptyset,\emptyset)(1)
 \\ & = & \int_{\alpha\in\suppfields {S^n}{\emptyset}{\emptyset}}
\KK {V^\dprime} (\alpha, \emptyset) \circ
\KK {V^\prime} (\emptyset,\alpha)(1)
\dsuppmeas {S^n} {\emptyset}{\emptyset} \\
 &=& \KK {V^\dprime} (0,\emptyset) \circ
\KK {V^\prime} (\emptyset, 0) (1) \;\;\;
\mbox{since $\suppfields {S^n} {\emptyset}{\emptyset} =
\{0\}$} \\
 &=& \KK {V^\dprime} (0,\emptyset) \circ \Id \circ
\KK {V^\prime} (\emptyset, 0) (1) \\
 &=& \KK {V^\dprime} (0,\emptyset) \circ
\KK {\op{D}} (\emptyset , 0 ) \circ \KK D (0,\emptyset)
\circ \KK {V^\prime} (\emptyset, 0) (1) \;\;\;
\mbox{by Lemma~\ref{lem:disks} } \\
 &=& Z^{[\theta]}_A(W^\prime)Z^{[\theta]}_A(W^\dprime).
\end{eqnarray*}
\end{proof}

\subsection{Invariants of products}

In this section we discuss the calculation of the invariants
of the product of two closed manifolds. Let $W$ and
$W^\prime$ be closed, oriented and connected of dimension
$m+1$ and $n+1$ respectively.

\begin{lem}
There is an identification of measure spaces
\[
\MF_{W\times W^\prime} \cong \MF_W \times \MF_{W^\prime}.
\]
\end{lem}

\begin{proof}
This follows from the fact that
$H_1(W\times W^\prime;\ZZ) \cong
H_1(W;\ZZ) \oplus H_1(W^\prime;\ZZ)$
and from the fact that all of these field spaces are given
the normalized Haar measure.
\end{proof}

Since $K_A$ is an $H$--space there is a Pontrjagin slant
product
\[
\slant\colon H^{m+n+2}(K_A;U(1)) \otimes H_{m+1}(K_A;\ZZ)
\ra H^{n+1}(K_A;U(1)).
\]
If $[W]\in  H_{m+1}(W;\ZZ)$ is the fundamental class then
given $\nu\in\MF_W$ we have $\nu_*[W] \in H_{m+1}(K_A;\ZZ)$.

\begin{thm}\label{thm:products}
Let $[\theta] \in H^{m+n+2}(K_A;U(1))$. Then
\[
Z^{[\theta]}_A (W\times W^\prime) = \int_{\nu\in \MF_W}
Z^{[\theta]\slant \nu_*[W]}(W^\prime) \dmeas W .
\]
\end{thm}

\begin{proof}
First recall that the slant product satisfies
\[
\langle a, b\bullet c \rangle = \langle a\slant b , c \rangle ,
\]
where $\bullet$ denotes the Pontrjagin product. Thus for
$v = (\nu,\nu^\prime) \in \MF_{W\times W^\prime} \cong
\MF_W \times \MF_{W^\prime}$
we have
\[
\langle v^*[\theta], [W]\times [W^\prime] \rangle =
\langle [\theta],
  \nu_*[W]\bullet \nu^\prime_* [W^\prime] \rangle =
\langle [\theta] \slant \nu_*[W],
  \nu^\prime_* [W^\prime] \rangle =
\langle \nu^{\prime *} ([\theta] \slant \nu_*[W]),
  [W^\prime] \rangle .
\]
Hence
\begin{eqnarray*}
Z^{[\theta]}_A (W\times W^\prime)  &=&
\int_{v\in \MF_{W\times W^\prime}} \langle v^*[\theta],
    [W]\times [W^\prime] \rangle \dmeas {W\times W^\prime} \\
 &=& \int_{\nu\in \MF_{W}}\int_{\nu^\prime
  \in \MF_{W^\prime}}
  \langle \nu^{\prime *} ([\theta] \slant \nu_*[W]),
    [W^\prime] \rangle \dmeas {W^\prime} \dmeas W \\
 &=& \int_{\nu\in \MF_W}
    Z^{[\theta]\slant \nu_*[W]}(W^\prime) \dmeas W .
\end{eqnarray*}
\end{proof}

\begin{exe}\label{ex:scfactor}
The product $M\times N$ where $M$ is simply connected.
{\em Let $M$ and $N$ be closed manifolds of dimension $m$
and $n$ respectively and let $\theta \in C^{m+n}(K_A;U(1))$
be a cocycle. For $0 \in H_m(K_A;\ZZ)$, we have
$[\theta] \slant 0$ trivial so
\[
 Z_A^{[\theta]} (M\times N) =
 \int_{\nu\in\MF_{M}}Z_A^{[\theta]\slant \nu_*[M]}(N)
 \dmeas {M} = Z_A^{[\theta]\slant 0}(N) = 1.
\]
}
\end{exe}

\begin{exe} \label{ex:smprod} The product $S^1\times M$.
{\em
Let $M$ be an $n$-manifold and let $\omega\colon A^{n+1} \ra U(1)$ be a
group cocycle corresponding to $\theta \in C^{n+1}(K_A; U(1))$. Noting
that $H_1(K_A;U(1))\cong A$ the slant product takes the form
\[
\slant \colon H^{n+1}(K_A;U(1))\otimes A \ra H^{n}(K_A;U(1))
\]
and may be described in terms of group cohomology as follows. For
$a\in A$ the slant product $\omega \slant a\colon A^n \ra U(1)$ is
given by
\begin{equation}\label{eq:slant}
(\omega \slant a)(g_1, \ldots , g_n) = \prod_{i=0}^n \omega(g_1,
\ldots , g_{i}, a, g_{i+1}, \ldots , g_n)^{(-1)^{\lambda_i}} .
\end{equation}
where $\lambda_i$ is the sign of the permutation taking
$(g_1, \ldots, g_n,a)$ to 
$(g_1,\ldots , g_{i}, a, g_{i+1}, \ldots , g_n)$.
This arises by using the
Eilenberg--Zilber map given by shuffle product.
Thus we have
\begin{equation}\label{eq:smprod}
 Z_A^{[\omega]} (S^1\times M) =
 \int_{a\in A}Z_A^{[\omega \slant a]}(M) d\!\mu
\end{equation}
and we can calculate an expression for the integrand using the
expression (\ref{eq:slant}) above.
}
\end{exe}


\section{Calculations}

Formulae such as that occurring in Theorems~\ref{thm:decomp}
and \ref{thm:products} are good tools for calculations. For
example we can fully compute all invariants in dimension 1+1
with almost no further effort. Suppose we have been given a
normalized group $2$--cocycle $\omega$ corresponding to the
defining cocycle $\theta \in C^2(K_A;U(1))$. We have already
computed the invariant for $S^2$. For $T^2$ we have
\begin{eqnarray*}
Z^{[\omega]}_A (T^2) = Z^{[\omega]}_A (S^1\times S^1)
 &=& \int_{a \in A} Z_A^{[\omega \slant a]}(M) d\!\mu \;\;\;
 \mbox{ by (\ref{eq:smprod})} \\
 &=& \int_{a \in A}\int_{b\in A} (\omega \slant a)(b)
        d\!\mu \;\;\; \mbox{ by Example \ref{ex:spheres}} \\
 &=& \int_{(a,b) \in A\times A} \omega(a,b)\op{\omega}(b,a)
        d\!\mu \;\;\; \mbox{ by (\ref{eq:slant})}.
\end{eqnarray*}
Finally, since a surface $\Sigma_g$ of genus $g$ is the
connected sum of $g$ tori we use
Proposition~\ref{prop:connectedsum} to get
\[
Z^{[\omega]}_A (\Sigma_g) =  Z^{[\omega]}_A (T^2)^g.
\]

One also needs to be able to make explicit calculations based on
explicit choices of the various cycles and cocycles in the
definitions. This takes us closer to the combinatorial view, but
it is important to remember that from the point of view of this
paper these are to be {\em deduced} not taken as definitions. This
is in fact the way Dijkgraaf and Witten introduced their
invariants: the path integral definition came first, followed by
the combinatorial formulae used to make explicit calculations.

\subsection{$\Delta$--complexes}

Everything in this section can be found elsewhere, but for
convenience we reproduce the essentials. It is convenient
for us to work with $\Delta$--complexes, as defined by
Hatcher \cite{Hatcher}, rather than simplicial complexes,
since $\Delta$--complexes will allow us to model manifolds
with far fewer simplices.

\begin{defn}
Suppose we have a collection of simplices $\{\Delta_i\}$,
together with an ordering (or numbering) of the vertices of
each simplex. As a result, we also get orderings on the sets
of vertices in the faces of the simplices $\Delta_i$. We can
now form a topological space by first taking the disjoint
union of the $\Delta_i$ and then identifying certain chosen
subsets $F_j$ of the faces of the $\Delta_i$ using the
canonical linear homeomorphisms that preserve the orderings
of the vertices (all faces in a given set $F_j$ are assumed
to be of the same dimension). A space which is constructed
in this way is called a {\em $\Delta$--complex.}
\end{defn}

Most of the ``triangulations'' of manifolds used in the
existing literature on DW--invariants are in fact
$\Delta$--complexes rather than simplicial complexes. The
same will apply in this paper, i.e.\ when we talk about a
triangulation of a manifold $M$, we mean a $\Delta$--complex
homeomorphic to $M$.  The main difference between a
$\Delta$--complex and a simplicial complex is that not every
simplex of a $\Delta$--complex has to be uniquely determined
by the set consisting of its vertices. The numbering of the
vertices in each simplex is needed to remove resulting
ambiguities.
Any simplicial complex can be turned into a $\Delta$--complex
by choosing an ordering of the vertices (this will induce an
ordering of the vertices of each simplex). Conversely, any
$\Delta$--complex is homeomorphic to a simplicial complex,
which can be constructed by subdivision of the simplices in
the $\Delta$--complex.

Homotopy classes of maps from a $\Delta$--complex $T$ to an
Eilenberg--Mac\,Lane space can be understood in combinatorial
terms as follows. A {\em colouring} of $T$ by the group $A$
is a map $g$ from the set of oriented edges of $T$ to $A$.
If $E$ is the oriented edge from the vertex labelled $a$ to
the vertex labelled $b$ (with $a<b$), then we denote $g(E)$
also as $g_{ab}$. We will use the convention that
$g_{ab}=g_{ba}^{-1}$ for all pairs of vertices $a,b$ which
are connected by an edge. Also, we impose a
{\em flatness condition}, which requires that, for any
triangle in $T$, the product of the colours on the boundary
is unity. More precisely, denoting the vertices of the triangle
by $a$, $b$ and $c$, we require that $g_{ab}g_{bc}=g_{ac}$, or
equivalently $g_{ab}g_{bc}g_{ca}=e$.
We define a {\em gauge transformation} to be a map $h$ from
the set of vertices of $T$ into $G$. We will often write
$h_a$ for $h(a)$. Gauge transformations form a group under
pointwise multiplication (in fact this group is isomorphic
to $G^{\MV}$, where $\MV$ is the number of vertices in $T$ ).
The group of gauge transformations has an action $\cdot$ on
the set of colourings, given by
\begin{equation}
(h\cdot g)_{ab}=h_b g_{ab}(h_a)^{-1}.
\end{equation}
The next proposition describes homotopy classes of maps from
a $\Delta$--complex to an Eilenberg--Mac\,Lane space
$K_A=K(A,1)$. Although it is well-known we include a proof
for completeness.

\begin{prop}\label{prop:colour}
Let $W$ be a manifold and $T$ a triangulation of $W$, then
the orbits of colourings of $T$ under gauge transformations
are in one to one correspondence with homotopy classes of
maps from $W$ into $K_A$. Moreover, homotopy classes of
{\em based} maps from $W$ to $K_A$ are in one to one
correspondence with orbits of colourings of $T$ under gauge
transformations which send a chosen vertex $x_0$ of $T$ to
the unit element of $A$.
\end{prop}

\begin{proof}
Start with a map $\sigma \colon W \ra K_A$. After a
suitable homotopy we can assume that $\sigma$ maps all
vertices of $T$ to the same point of $K_A$. Hence all
edges of $T$ become loops in $K_A$, and since
$\pi_1(K_A)\cong A$ we can color each edge of $T$ with an
element of $A$. All colourings of $T$ induced in this way
satisfy a flatness condition because the image of any
triangle in $K_A$, and hence also the image of the loop which
forms its boundary, is contractible. One should note that one
may obtain different colourings of $T$ from the same homotopy
class of maps. It is easy to see why this happens. Suppose
that we have two homotopic maps $\sigma$ and $\sigma'$ from
$W$ to $K_A$ which both send all vertices of $T$ to the base
point for $\pi_1(K_A)$. Although $\sigma$ and $\sigma'$ are
homotopic, the homotopy between them may move the vertices of
$T$ around non-contractible loops in $K_A$. If the vertex $v$
gets moved around the loop labelled by $h \in A$, then the
group elements of the edges of $T$ which end at $v$ get
multiplied by $h$ from the left, while the group elements on
edges which begin at $v$ get multiplied by $h^{-1}$ from the
right. This is exactly the effect of a gauge transformation
at the vertex $v$. Thus we do not get a well-defined map from
homotopy classes of maps to colourings of $T$, but we do get
a well-defined map from homotopy classes of maps to gauge
orbits of colourings of $T$. This map is in fact invertible.
To see injectivity, suppose that two maps $\sigma$ and
$\sigma'$ induce the same gauge class of colourings of $T$.
Then these maps are certainly homotopic on the $1$--skeleton
of $T$ and, using the fact that $K_A$ has trivial higher
homotopy, we may extend the homotopy on the $1$--skeleton to
a homotopy on all of $T$, or $W$. For surjectivity, take
any colouring of $T$ satisfying the flatness condition. We
may always construct a map from the $1$--skeleton of $T$
into $K_A$ which induces this colouring and, because $K_A$
has trivial higher homotopy, this map extends to a map from
all of $W$ to $K_A$. The statement about based maps follows
in a similar way if we identify the base point of $W$ with
the chosen vertex $x_0$ of $T$. This vertex can now no
longer be moved around $K_A$ by homotopies and hence
colourings, which differ by a non trivial gauge
transformation at $x_0$, do not correspond to the same
homotopy class of based maps.
\end{proof}

\subsection{A formula for explicit calculation}

We will assume that the HQFT defining the abelian homotopy
Dijkgraaf--Witten theory is integrable (see
Sect.~\ref{subsec:measure}) and that we have a group cocycle
$\omega$ corresponding to the defining singular cocycle
$\theta$. Given $\alpha_0\in \MF_{M_0}$ and
$\alpha_1 \in \MF_{M_1}$ we need to determine the effect of
the maps
$E_{W,\nu} \colon L_{M_0,\alpha_0} \ra L_{M_1,\alpha_1}$
occurring in (\ref{eq:matrix}), where
$\nu\in \MFr W {\alpha_0}{\alpha_1}$.
To do this we must make some choices:
\begin{itemize}
\item choose a representative $\sigma\colon W \ra K_A$
      of the class $\nu$,
\item choose fundamental cycles $a_i\in C_nM_i$ for $i=0,1$
      (giving generators of $L_{M_i,\gamma_i}$, where
      $\gamma_i = \sigma|_{M_i}$),
\item choose $f\in C_{n+1}W$ representing the fundamental
      class in $H_{n+1}(W,\partial W)$.
\end{itemize}
Armed with these choices we then compute $\sigma^*\theta (f)$
and hence determine
$E_{W,\sigma}(a_0) = \sigma^*\theta (f)a_1$.

Let us now suppose that T is a triangulation of $W$ which
induces triangulations $T_0$ and $T_1$ of $M_0$ and $M_1$.
Since $T$, $T_0$ and $T_1$ are $\Delta$--complexes, they
immediately give canonical representatives $f$, $a_0$ and
$a_1$ for the fundamental classes of $W$, $M_0$ and $M_1$
and moreover these satisfy $\partial f = a_1 - a_0$.
Explicitly, for $i=0,1$ we have
\begin{equation}
\label{fundtet} a_i=\sum_{t \in T_i} \epsilon_t [t],
\end{equation}
where the sum runs over the $n$--simplices of $T_i$ and
$[t]$ denotes the inclusion map of the $n$--simplex $t$ into
$T_i$ (i.e.\ the inclusion map into the set of disjoint
simplices followed by the identification map). The signs
$\epsilon_t$ express the orientation of the simplices
compared to that of the whole manifold. Note that the
orientation of a simplex can be described in terms of the
ordering of its vertices. Hence the signs $\epsilon_t$ are
fixed by the orientation of $M$ and the chosen
$\Delta$--complex structure. Similarly we have
\begin{equation}
f=\sum_{t \in T} \epsilon_t [t],
\end{equation}
where here the sum is over the $n+1$--simplices of $T$.

Next, given $\nu \in \MF_W$ we use
Proposition~\ref{prop:colour} to choose a colouring of $T$
(in general there may be many such colourings). Now define
a map $\sigma \colon W \ra K_A$ such that $[\sigma] = \nu$
as follows.
Choose representatives for the elements of the fundamental
group of $K_A$, or more precisely, for every $g \in A$ fix
a map $l_g$ from the standard $1$--simplex onto a loop in
$K_A$ which corresponds to the element
$g \in \pi_1(K_A)\cong A$. Using these, we can define
$\sigma$ on the $1$--skeleton of the triangulation by
mapping an edge labelled $g$ into $K_A$ by $l_g$.  To fix
$\sigma$ on the $2$--skeleton one introduces standard maps
from any coloured $2$--simplex to $K_A$, such that these
maps reduce to the standard maps for $1$--simplices on
the coloured boundary. One continues in this way for the
higher skeleta until $\sigma$ is defined (these map
extensions are possible because $K_A$ has trivial higher
homotopy). It is clear by the proof of
Proposition~\ref{prop:colour} that $[\sigma] = \nu$.

If $t$ is an $n+1$--simplex in $T$ then $\sigma^*\theta (t)$
is a function of the colouring chosen above and we can assume
that $\theta$ and $\omega$ are related so that
\[
\sigma^*\theta (t) =
\omega(g^{\sigma}_{t,1},\ldots,g^{\sigma}_{t,n+1}),
\]
where $g^\sigma_{t,1},\ldots,g^\sigma_{t,n+1}$ are the group
elements which colour $n+1$ edges which don't lie in the same
face (flatness then determines the others). We will take
these $n+1$ edges to be the edges which connect the vertices
of the simplex in ascending order\footnote{Note that if we
were only given $\sigma$, the procedure described here gives a
way of determining a suitable $\omega$.}.
Thus (using multiplicative notation for the group operation
in $U(1)$) we have
\[
\sigma^*\theta(f) = \sigma^*\theta (\sum_{t \in T}
\epsilon_t [t])
= \prod_{t\in T} \sigma^*\theta (t)^{\epsilon_t}
= \prod_{t\in T}\omega(g^{\sigma}_{t,1},
       \ldots,g^{\sigma}_{t,n+1})^{\epsilon_t}.
\]
When $W$ is closed, the number $\sigma^*\theta(f)$ does not
depend on the chosen triangulation of $W$ (which corresponds
to a choice of $f$) or on the choice of $g^{\sigma}$ in its
gauge orbit. If $W$ is not closed, then we will still have
the same formula as above, but, since $f$ has non-zero
boundary in this case, the number $\sigma^*\theta(f)$ will
now depend on the choice of $\sigma$, as well as on the
choice of $f$, that is, of the triangulation. Nevertheless,
one may check that any choice would still determine the same
map $E_{W,\nu}$.

If we choose the same $a_0$, $a_1$ and $f$ for each
$\nu \in \MFr W {\alpha_0}{\alpha_1}$ then
$\KK W (\alpha_0,\alpha_1) $ is described by
\begin{equation}\label{eq:invstatesum}
\KK W (\alpha_0,\alpha_1) (a_0) =
(\int_{\nu = [\sigma]\in \MFr W {\alpha_0}{\alpha_1}}
\prod_{t\in T}
\omega(g^{\sigma}_{t,1},\ldots,g^{\sigma}_{t,n+1})^{\epsilon_t}
    \dmeasalpha ) ~a_1.
\end{equation}
For a closed $n+1$--manifold we get
\begin{equation}\label{eq:triazfrm}
Z_{A}^{[\theta]}(W)=\int_{\nu =
 [\sigma]\in\MF_W} \prod_{t\in T}
\omega(g^{\sigma}_{t,1},\ldots,g^{\sigma}_{t,n+1})^{\epsilon_t}
\dmeas W .
\end{equation}
Note that there is nothing in the above depending on any
special property of the group $A$. As long as a good measure
on the space of homotopy classes of based maps
$\MF_W=[W;K(A,1)]$ is available, the above formulae can be
used to calculate the invariants. The reason for restricting
to compact abelian Lie groups $A$ is that we have good
measures available as already stated in the introduction.
Of course for finite $A$ one also has a measure (the
counting measure) available in case $A$ is not abelian,
and in the state sum approach one actually starts with
the above formulae (\ref{eq:invstatesum}) and
(\ref{eq:triazfrm}) for the invariants.

\subsection{Dimension 2+1}

In this last section we take $A=U(1)$ and at level $k$ we
use the group cocycle $\omega_k$ defined in
(\ref{eq:gpcocycle}). We will write $Z^k(W)$ to mean
$Z^{[\omega_k]}_{U(1)}(W)$ and by the ``$U(1)$ homotopy
DW--invariants'' of a closed $3$--manifold
$W$ we mean the collection of numerical invariants
$\{ Z^k(W)\}_{k\geq 0}$. We will prove the following
proposition.

\begin{thm}\label{thm:lens}
The $U(1)$ homotopy Dijkgraaf--Witten invariants
distinguish homotopy equivalence classes of lens spaces.
\end{thm}

Before proving this let us recall certain facts about lens
spaces. Lens spaces are a class of $3$--manifolds
parametrized by pairs of coprime integers $(p,q)$, the lens
space labelled by $(p,q)$ being denoted $L(p,q)$. 
Since we are interested here in
oriented and not only orientable lens spaces a bit of care
is needed. Our orientation convention will be the standard
one, i.e.\ $L(p,q)$ is the closed oriented $3$--manifold
obtained by surgery on $S^{3}$ along the unknot with surgery
coefficient $-p/q$, where $L(p,q)$ is given the orientation
induced by the standard right-handed orientation on $S^{3}$.
We note that
\begin{itemize}
\item The lens spaces $L(p,q)$ and $L(p',q')$ are
      homeomorphic if and only if $p$ is equal to $p'$ and
      $q=\pm q' \bmod p$ or $qq'=\pm 1 \bmod p$.
\item $L(p,q)$ and $L(p',q')$ are homotopy equivalent if and
      only if $p=p'$ and $qq'=\pm a^2 \bmod p$ for some
      integer $a$.
\end{itemize}
The first fact was proved by Reidemeister, cf.\
\cite{Reidemeister}, and the second fact is due to Whitehead
\cite{Whitehead}. For a more recent source, see for instance
\cite{Rolfsen, Saveliev}. In all cases, the minus sign
corresponds to a reversal of the orientation.  We will be
interested in homotopy classes of lens spaces using only
orientation preserving homeomorphisms, since the DW--invariants
depend on the orientation (e.g.\ they can
have different values for, say, $L(p,q)$ and $L(p,p-q)$).
Therefore, in the rest of the paper, when we say that two
lens spaces $L(p,q)$ and $L(p,q')$ are homotopy equivalent,
this means that $qq'=+a^2 \bmod p$ for some $a \in \ZZ$. We
note that $L(0,\pm 1)=S^{2} \times S^{1}$ with fundamental
group $\ZZ$. All the abelian homotopy DW--invariants
of this manifold are trivial by (\ref{eq:smprod})
and Ex.~\ref{ex:spheres} (alternatively use
Ex.~\ref{ex:scfactor}). From now on we assume that $p \neq 0$.
Note then that the fundamental group of $L(p,q)$ is $\ZZ/p$
and the other homotopy groups are isomorphic to those of the
$3$--sphere. Hence the homotopy groups of a lens space do not
determine its homotopy type.

\begin{figure}[h,t,b]
\centerline{ \psfig{figure=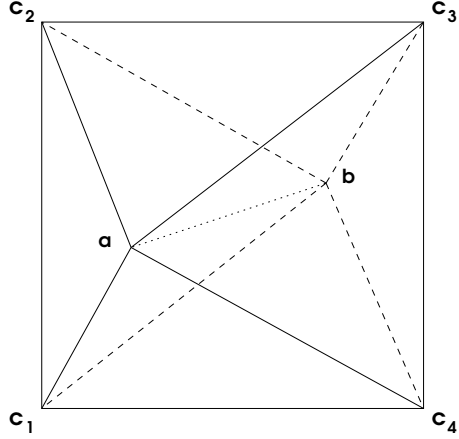,width=6cm}}
\caption[something]{\footnotesize The polyhedron from which
$L_{4,1}$ is formed by identification of each face on the
front with the next face on the back} \label{lensfig}
\end{figure}

The lens space $L(p,q)$ has a nice triangulation consisting
of $p$ tetrahedra with vertices
$a_i$, $b_i$, $c_i$ and $d_i$, $i=1,\ldots,p$, illustrated for $p=4$
in Figure~\ref{lensfig}. The tetrahedra are first glued
together along the $abc$--faces, i.e.\ we make the
identification
$(a_i, b_i, c_i) \equiv (a_{i+1},b_{i+1},c_{i+1})$ for all
$i$ with the convention that $a_{p+1}=a_{1}$ etc.

After these identifications there is one point corresponding
to all the $a_i$, which we will call $a$ and there is
similarly one point corresponding to the $b_i$ denoted $b$.
To get the lens space $L_{p,q}$ from this polyhedron, one
identifies each face on one side with the face which lies
$q$ steps clockwise removed on the other side, i.e.\
one makes the identification
$(a,c_{i}, d_{i})\equiv (b,c_{i+q},d_{i+q})$, again with
$c_{p+1}=c_1$ etc. The path $ab$ has now become a loop
and one may easily check that it is a generator of the
fundamental group. One may number the vertices such that
the signs $\epsilon_t$, which occur in the formula for the
fundamental cycle, are all positive.

For the U(1) homotopy Dijkgraaf--Witten theory 
the space of fields is
\[
\MFLpq = H^1(L(p,q);U(1)) = \Hom (\ZZ/p, U(1)) =
\{\zeta \in U(1) \; | \; \zeta^p = 1\} =: \Lambda_p .
\]
Colourings of the above triangulation were studied by
Altschuler and Coste \cite{AltschulerCoste}
(for finite groups which is sufficient here as
$\Lambda_p \cong \ZZ/p$). Given
$\nu \in \Lambda_p $ they provide a particularly nice
colouring corresponding to $\nu$ by colouring
the three independent edges (in ascending order) in the $j$'th
tetrahedron $t_j$ with the group elements $\nu, \nu^{j\qb}$
and $\nu^{\qb}$ respectively, where $\qb$ is the inverse of
$q$ modulo $p$. Using (\ref{eq:triazfrm}) we then have
\[
Z^k(L(p,q))=\int_{\nu \in\MFLpq} \prod_{j=1}^p
\omega_k(\nu, \nu^{j\qb}, \nu^{\qb}) \dmeas L
= \frac{1}p \sum_{\nu\in \Lambda_p}
  \prod_{j=1}^p  \omega_k(\nu, \nu^{j\qb}, \nu^{\qb}).
\]
For $u\in U(1)$ let $\langle u \rangle$ be the unique number
in the interval $[0,1)$ such that
$u = e^{2\pi i \langle u \rangle}$. It is easy to see that
$\sum_{j=1}^p \langle \nu^{j\qb} \rangle =
\sum_{j=1}^p \langle \nu^{\qb (j+1)} \rangle $
and so we can write
\begin{eqnarray*}
\prod_{j=1}^p \omega_k(\nu, \nu^{j\qb}, \nu^{\qb}) &=&
\prod_{j=1}^p e^{2 \pi i k \langle \nu \rangle (\langle
  \nu^{j\qb}\rangle+ \langle \nu^{\qb} \rangle
       - \langle \nu^{\qb (j+1)}\rangle)} \\
 &=& e^{2 \pi i k \langle \nu \rangle \sum_{j=1}^p
   (\langle \nu^{j\qb}\rangle+ \langle \nu^{\qb} \rangle
    - \langle \nu^{\qb (j+1)}\rangle)} \\
 &=& e^{2 \pi i k \langle \nu \rangle p
    \langle \nu^{\qb} \rangle}\\
 &=& e^{\frac{2 \pi i k  \qb l^2}{p}},
\end{eqnarray*}
where in the last equality we have written
$\nu = e^{\frac{2 \pi i l}{p}}$ for some $l=1 \ldots p$.
Thus we have
\begin{equation}\label{eq:lenscalc}
Z^k(L(p,q))=\frac{1}p \sum_{l=1}^{p}
      e^{\frac{2 \pi i k  \qb l^2}{p}}.
\end{equation}
Let us recall formulas for the involved Gauss sums. For
$r,N$ relatively prime, let us write
\begin{equation}
G(r,N):= \sum_{l=1}^{N} e^{\frac{2 \pi i r l^2}{N}}.
\end{equation}
Dirichlet \cite{Dirichlet1,Dirichlet2} proved that
\begin{equation}
\label{dirfrm}
G(r=1,N)=\left\{ \begin{array}{ll}
(1+i)\sqrt{N}, & N=0 \bmod 4, \\
\sqrt{N}, & N=1 \bmod 4,\\
0, & N=2 \bmod 4,\\
i \sqrt{N}, & N=3 \bmod 4.
\end{array}\right.
\end{equation}
Futhermore, when $N$ is an odd prime, there is a closed
formula for $G(r,N)$ for all $r$,
\begin{equation}
\label{primefrm}
G(r,N)=\left\{ \begin{array}{ll}
\leg{r}{N}\sqrt{N}, & N=1 \bmod 4, \\
i \leg{r}{N}\sqrt{N}, & N=3 \bmod 4,
\end{array}\right.
\end{equation}
where $\leg{r}{N}$ is the Legendre symbol or $r$ modulo $N$,
that is, $\leg{r}{N}$ equals $1$ if $r$ is a square modulo $N$
and $-1$ otherwise.
Before proving Theorem~\ref{thm:lens} we require the
following lemma.

\begin{lem}\label{classprop}
Let $p=2^{k}p_{1}^{k_1}p_{2}^{k_2}\ldots p_{m}^{k_m}$ be the
prime decomposition of $p$ {\em (}the $p_i$ are odd primes,
$k$ is nonnegative and the $k_{i}$ are positive{\em)} and
consider homotopy classes of lens spaces $L(p,q)$. We
distinguish three cases.
\begin{itemize}
\item{\em $k=0$ or $k=1$.} There are $2^m$ homotopy classes
which we can label by the string of signs
$\left(\leg{q}{p_1},\ldots,\leg{q}{p_m}\right)$.
\item{\em $k=2$.} There are $2^{m+1}$ homotopy classes which
may be labelled by $q \bmod 4$ and the signs $\leg{q}{p_i}$.
{\em (}Note that $q \bmod 4$ equals $1$ or $3$.{\em )}
\item{\em $k>2$.}  There are $2^{m+2}$ homotopy classes
labelled by $q \bmod 8$ and the signs $\leg{q}{p_i}$.
{\em (}Note that $q \bmod 4$ equals $1$, $3$, $5$ or $7$.
{\em )}
\end{itemize}
\end{lem}

\begin{proof}
Recall that $\ZZ^{*}_{p}$, the multiplication group
modulo $p$, decomposes as
\begin{equation}\label{zpdec}
\ZZ^{*}_p = \ZZ^{*}_{2^{k}}\times\ZZ^{*}_{p_{1}^{k_1}}
\times\ldots\times\ZZ^{*}_{p_{m}^{k_m}}.
\end{equation}
Hence, if $x$ is an element of $\ZZ_{p}^{*}$ we may
write $x = (x_{0},x_{1},\ldots,x_{m})$ with
$x_{i} \in \ZZ_{p_{i}^{k_i}}^{*}$ (with $p_0 =2$, $k_0 =k$).
In fact, we can take
$x_{i}=x \bmod p_{i}^{k_i}$. From this decomposition it is
clear that $x$ will be a square modulo $p$ if and only if
$x$ is a square modulo $p_{i}^{k_i}$ for
$i=0,1,\ldots,m$. Furthermore, it is not difficult to show that
$x$ is a square modulo $2^k$ if and only if $x=1 \bmod 8$ and
$x$ is a square modulo $p_{i}^{k_i}$ if and only if $x$ is a
square modulo $p_i$, $i=1,\ldots,m$. 
To find the homotopy classes of lens
spaces we must therefore find out which elements of 
$\ZZ_n^{*}$ give a square when they are multiplied together,
$n$ being any odd prime. Obviously the
product of two squares is always a square. Also, using the
fact that $\ZZ^{*}_{n}$ is cyclic, one
sees that the product of two non-squares is a square in the
$\ZZ^{*}_{n}$, while the product of a square and a non-square
in $\ZZ^{*}_{n}$ is never a square. Finally we note that
two elements multiply to a square in $\ZZ^{*}_{2^k}$
only if they are equal modulo the
minimum of $8$ and $2^k$.
\end{proof}

\begin{proof} (of Theorem~\ref{thm:lens})\\
For any lens space $L(p,q)$ we have from (\ref{eq:lenscalc})
that $Z^0(L(p,q)) = 1$. The next value of $k$ for which
$Z^k(L(p,q)) = 1$ occurs when $k=p$ (essentially this is the
triangle inequality for complex numbers), so this determines
$p$.

Now fix $p$ and write its prime decomposition as in
Lemma~\ref{classprop}. We need to show that the invariants
$Z^{k}$ determine the labels of the homotopy classes given
in that lemma. Let $p_i$ be one of the odd prime factors
(if there are no odd prime factors, we only need to determine
$q \bmod 4$ or $q \bmod 8$, see further on for that) and
consider $k=\frac{p}{p_i}$. Filling in (\ref{eq:lenscalc}),
we get
\begin{equation}
Z^{p/p_i}(L(p,q))= \frac{1}{p}\sum_{l=1}^{p}
\exp\left(\frac{2 \pi i\bar{q}l^2}{p_i}\right)
 =\frac{1}{p_i}\sum_{l=1}^{p_i}
\exp\left(\frac{2 \pi i\bar{q}l^2}{p_i}\right)
\end{equation}
and using (\ref{primefrm}), we see that
\begin{equation}
Z^{p/p_i}(L(p,q))=\left\{ \begin{array}{ll}
\frac{1}{\sqrt{p_{i}}}\leg{q}{p_{i}}, & p_{i}=1 \bmod 4, \\
\frac{i}{\sqrt{p_{i}}}\leg{q}{p_{i}}, & p_{i}=3 \bmod 4.
\end{array}\right.
\end{equation}
Thus these invariants determine the Legendre symbols
$\leg{q}{p_{i}}$. This means they separate homotopy classes
of lens spaces with $p$ odd or $p=2 \bmod 4$, the first case
in Lemma~\ref{classprop}. To settle the second case
($p=4 \bmod 8$), we need to determine $q \bmod 4$. This
is accomplished by taking $k = p/4$. We have
\begin{equation}
Z^{p/4}(L(p,q))=\frac{1}{4}\sum_{l=1}^{4}
e^\frac{2 \pi i\bar{q}l^2}{4}=\frac{1}{2}(1+i^{\bar{q}})=
\left\{ \begin{array}{ll}
\frac{1}{2}(1+i), & q=1 \bmod 4, \\
\frac{1}{2}(1-i), & q=3 \bmod 4.
\end{array}\right.
\end{equation}
To deal with the final case ($p=0 \bmod 8$), we have to
determine $q \bmod 8$. This can be done using $k ={p/8}$:
\begin{equation}
Z^{p/8}(L(p,q))=\frac{1}{8}\sum_{l=1}^{8}
e^\frac{2 \pi i\bar{q}a^2}{8}=
\frac{1}{4}(1+(-1)^{\bar{q}}+2e^\frac{\pi i\bar{q}a^2}{4})=
\left\{ \begin{array}{ll}
\frac{1}{2}e^{i \pi/4}, & q=1 \bmod 8, \\
\frac{1}{2}e^{3i \pi/4}, & q=3 \bmod 8, \\
\frac{1}{2}e^{5i \pi/4}, & q=5 \bmod 8, \\
\frac{1}{2}e^{7i \pi/4} & q=7 \bmod 8.
\end{array}\right.
\end{equation}
\end{proof}


\noindent {\bf Acknowledgements} All authors were supported by the
European Commission, JKS and PRT through Marie Curie fellowships
and SKH through the research network EDGE. SKH thanks the
School of Mathematics at the University of Edinburgh, and the
Max--Planck--Institut f\"{u}r Mathematik in Bonn for their
hospitality. He also was supported by the Danish Natural Science 
Research Council and the Max--Planck--Institut. JKS thanks the
School of Mathematical and Computer Sciences at Heriot--Watt
University and PRT thanks the Institut de Recherche 
Math\'ematique Avanc\'ee in Strasbourg.


\end{document}